\documentclass[twocolumn]{revtex4}
\usepackage{latexsym}
\usepackage{amsmath}
\usepackage{times}
\usepackage{amssymb}
\usepackage[T1]{fontenc}
\usepackage[utf8]{inputenc}
\usepackage{fancyhdr}
\usepackage{url}
\usepackage{hyperref}
\usepackage{color}
\usepackage{graphicx}

\usepackage{bbm}
\newcommand{\R}{\mathbb R}

\newcommand{\ti}{t_{\textup{i}}}
\newcommand{\tf}{t_{\textup{f}}}
\newcommand{\Dt}{\Delta{t}}
\newcommand{\h}{\Delta{x}}
\graphicspath{{figures/}}

\begin{document}

\title{\bf Blending Brownian Motion and Heat Equation}
\author{Emiliano Cristiani$^{1}$}\thanks{Corresponding author. Email: e.cristiani@iac.cnr.it}
\affiliation{$^1$Istituto per le Applicazioni del Calcolo, Consiglio Nazionale delle Ricerche, \\ Via dei Taurini 19, 00185 Rome, Italy}
\date{\today}


\begin{abstract}
\noindent In this short communication we present an original way to couple the Brownian motion and the heat equation. More in general, we suggest a way for coupling the Langevin equation for a particle, which describes a single realization of its trajectory, with the associated Fokker-Planck equation, which instead describes the evolution of the particle's probability density function. 
Numerical results show that it is indeed possible to obtain a regularized Brownian motion and a Brownianized heat equation still preserving the global statistical properties of the solutions. The results also suggest that the more macroscale leads the dynamics the more one can reduce the microscopic degrees of freedom.

\vspace*{2ex}\noindent\textit{\bf Keywords}: Brownian motion, heat equation, diffusion equation, multiscale methods, statistical properties, coupling.
\\[3pt]
\noindent\textit{\bf MSC}: 35K05, 60J65, 35Q84, 35Q40
\end{abstract}

\maketitle

\thispagestyle{fancy}

\section{Introduction}\label{sec:intro}
\noindent In this short communication we present an original way to couple the Brownian motion and the heat equation. 
The method can be used to describe, in a new multiscale fa\-shion, diffusion-based multiscale phenomena and dual-nature phenomena, like, e.g., the wave-particle duality of elementary particles in quantum mechanics.

This paper falls in the context of multiscale methods for the numerical solution of problems based on (stochastic) ODEs and PDEs arising in mathematics, physics, and engineering. 
Relevant literature is huge, see, e.g., the book \cite{e2011book} for a quick overview.
Although details can change considerably from method to method, a general idea behind multiscale modeling is to provide different levels of mathema\-tical description for phenomena occurring at different scales. 
In 2011 the paper \cite{cristiani2011MMS} proposed a modeling technique for advection problems whereby the micro and macro levels are advanced in time concurrently in the whole domain, and they are coupled together via a ``blending'' parameter, say $\theta\in[0,1]$, which measures the relative weight of one description over the other, say $\theta=0$ for fully macro and $\theta=1$ for fully micro. In different words, one defines a mesoscopic blend of two replicas of the same system, similar in spirit, although not in mathematical content, to mesoscopic approaches to multiscale modeling. 
The main idea behind this approach is to keep the two replicas in constant interaction. This allows to describe systems where the microscopic level must be observed everywhere and at any time, because, e.g., granular properties are never and nowhere negligible. This kind of systems can be found in several contexts like crowd dynamics, financial mathematics, and complex systems in general (consider, e.g. the well-known \textit{herding behavior}, where the choice of single individuals influence the mass and vice versa). The price to pay is the computational cost of keeping the two solvers active always and everywhere, a cost that can hopefully be optimized by minimizing the number of microscopic degrees of freedom.
Besides the computational aspects, however, this ``replica'' approach has the conceptual appeal of treating the micro and macro levels at the same footing, which might prove advantageous to describe systems which display genuine physical duality, like the wave-particle duality in quantum mechanics.
\section{Coupling Brownian motion and heat equation}
\noindent In this section we derive a way to couple the Brownian motion and the heat equation. 
For the sake of clarity, in the following we restrict the discussion to dimension one. 

\subsection{Brownian motion and heat equation}\label{sec:B&H}
\noindent Let us consider the following system of stochastic differential equations 
\begin{equation}\label{SDEs}
\left\{
\begin{array}{ll}
dX^k_t=\sqrt{2D} dW^k_t \\
X^k_0=0
\end{array}
\right.
\qquad k=1,\,\dots,\,N_p,
\end{equation}
where $D>0$ is the diffusion coefficient, $W^k_t$ is the one--dimensional Brownian motion (standard Wiener process), $N_p$ is the number of particles under observation, and $\{X^k_t\}_k$ their positions at time $t$. The associated Fokker--Planck equation is the heat equation
\begin{equation}\label{heateq_parab}
\left\{
\begin{array}{ll}
\frac{\partial}{\partial t} u-D\frac{\partial}{\partial x^2}u=0, & t>0,\ x\in\R \\[2mm]
u(0,\,x)=\delta_0, & x\in\R,
\end{array}
\right.
\end{equation}
where $\delta_0$ is the Dirac delta centered in 0. The solution $u$ to (\ref{heateq_parab}) describes the time evolution of the probability density function of the Brownian motions generated by (\ref{SDEs}) \cite{albert}, i.e.\
$$
\mathbb P(X^k_t\in \Omega)=\int_\Omega u(t,x)\,dx,
\quad \forall\,\Omega\subseteq\R,
\ \forall k.
$$ 
The solution to the heat equation (\ref{heateq_parab}) is 
\begin{equation}\label{eq:sol_esatta_calore_1d} 
u(t,x)=\frac{1}{\sqrt{4\pi D t}}e^{-\frac{x^2}{4D t}},\qquad t>0,\ x\in\R.
\end{equation}
For our purposes, it is useful to define a \emph{macroscopic velocity field} $v_M$ as
\begin{equation}\label{def_vM}
	v_M(t,x):=-D\frac{\frac{\partial}{\partial x}u(t,x)}{u(t,x)},
\end{equation}
so that the heat equation can be formally written as
\begin{equation}\label{heateq_iperb}
\frac{\partial}{\partial t}u+\frac{\partial}{\partial x}(u v_M)=0.
\end{equation}
The form (\ref{heateq_iperb}) makes the velocity field $v_M$, which transports the probability density function, appear explicitly. Note that this ``hyperbolic'' formulation is valid if $u\neq 0$, i.e.\ if $t>0$.
By (\ref{eq:sol_esatta_calore_1d}) and (\ref{def_vM}) we have 
\begin{equation}\label{Vm_esplicito}
v_M(t,x)=\frac{x}{2t}.
\end{equation}


%
%
%
%
%

\subsection{Numerical approximation}\label{sec:numericalapproximation}
\noindent Let us introduce a structured space-time grid with spatial cell size $\h$ and a time step $\Dt$. We denote the grid nodes by $(t_n,x^j)$, for $n=1,\ldots,N_t$ and $j=1,\ldots,N_x$, the space cells by $E^j=\left[x^j-\frac{\h}{2},\,x^j+\frac{\h}{2}\right)$, and by $f^j_n$ the approximate value of a generic function $f(t,x)$ at $t=t_n$ and $x=x^j$. 

To avoid to manage the initial datum $\delta_0$ at the discrete level, we set the initial time for the simulation at $t=\ti>0$ and, consequently, the new initial datum $u(\ti,x)=\frac{1}{\sqrt{4\pi D \ti}}\exp\left(-\frac{x^2}{4D \ti}\right)$. In this way we can safely assume that $u>0$ for any $t\geq \ti$.
The computational domain is set at $[\ti,\tf]\times[-x_\textup{b},x_\textup{b}]$, with $\ti=\frac12$, $\tf=5$, and $x_\textup{b}=8$. We choose $N_x=50$ ($\Delta x=0.32$) and $N_t=500$ ($\Delta t=0.009$). Without loss of generality, in the numerical tests we always set $D=\frac12$.

\medskip

The Ito processes satisfying (\ref{SDEs}) can be approximated by using the weak Euler scheme
\begin{equation}\label{Euler_for_SDEs}
	Y^k_{n+1}=Y^k_n+\omega_n\sqrt{2D}\sqrt{\Delta t}, \qquad k=1,\,\dots,\,N_p
\end{equation}
where $\{\omega_n\}_n$ are independent two-point random variables with $\mathbb P(\omega_n=\pm 1)=\frac{1}{2}\ \forall n$.
The discrete velocities of the particles can be computed \emph{a posteriori} as $\frac{\omega_{n}\sqrt{2D}\sqrt{\Delta t}}{\Dt}$ (note that the velocity tends correctly to $+\infty$ as $\Delta t \to 0$). 
At time $t=\ti$, we assume that microscopic particles $\{Y^k_{\ti}\}_k=\{X^k_{\ti}\}_k$ are distributed with probability density $u(\ti,x)$.

Given the particles' positions, we recover the approximate probability density function by
\begin{equation}\label{defpsi}
\psi_n^j:=\frac{1}{N_p \Delta x}\sum_{k=1}^{N_p} \mathbbm{1}_{E^j}(Y^k_n), \quad \forall n, \forall j,
\end{equation}
where $\mathbbm{1}$ denotes the characteristic function. Note that $\frac{1}{N_p}$ is the mass carried by each particle, the total mass being $\int_\R u dx=1$.

\medskip

Equation (\ref{heateq_parab}) can be easily discretized by means of a first-order centered-in-space and forward-in-time finite difference scheme
\begin{equation}\label{schemeDFC}
	\phi_{n+1}^{j}=\phi_{n}^{j}+D\frac{\Dt}{\h^2}
		\left(\phi_{n}^{j+1}-2\phi_{n}^{j}+\phi_{n}^{j-1}\right).
\end{equation}

\subsection{Partial coupling I: The regularized Brownian motion}\label{sec:partialcouplingmicro}
\noindent A preliminary coupling of the microscopic and macroscopic scale can be obtained by assuming that the solution $u$ to the heat equation (\ref{heateq_parab}) is known, see (\ref{eq:sol_esatta_calore_1d}), and then solving the equation for $X^k$'s, duly coupling the original microscopic velocity field and the macroscopic velocity field (\ref{Vm_esplicito}). Following slavishly \cite{cristiani2011MMS}, we should blend the two velocity fields by a convex combination getting
\begin{equation}\label{eq:partialcouplingmicroNO}
	Y^k_{n+1}=Y^k_n+\theta\omega_{n}\sqrt{2D}\sqrt{\Delta t}+(1-\theta)\frac{Y^k_n}{2t_n}\Dt,
\end{equation}
for $k=1,\,\dots,\,N_p$, $\theta\in[0,1]$. Surprisingly enough (or not?), this choice does not give the desired result, in the sense that the corresponding probability density function $\psi$ does not approximate correctly the function $u(t,x)$ for $\theta\in(0,1)$. A good scale interpolation is obtained instead by setting, at any time step $n$ and for any particle $k$,
\begin{equation}\label{eq:partialcouplingmicro}
Y^k_{n+1}=
\left\{
\begin{array}{ll}
Y^k_n+\omega_{n}\sqrt{2D}\sqrt{\Delta t}, & \textrm{with prob.\ } \theta, \\ [2mm]
Y^k_n+\frac{Y^k_n}{2t_n}\Dt, & \textrm{with prob.\ } (1-\theta).
\end{array}
\right.
\end{equation}

\begin{figure}[h!]
\begin{center}
\begin{tabular}{cc}
\includegraphics[width=0.23\textwidth]{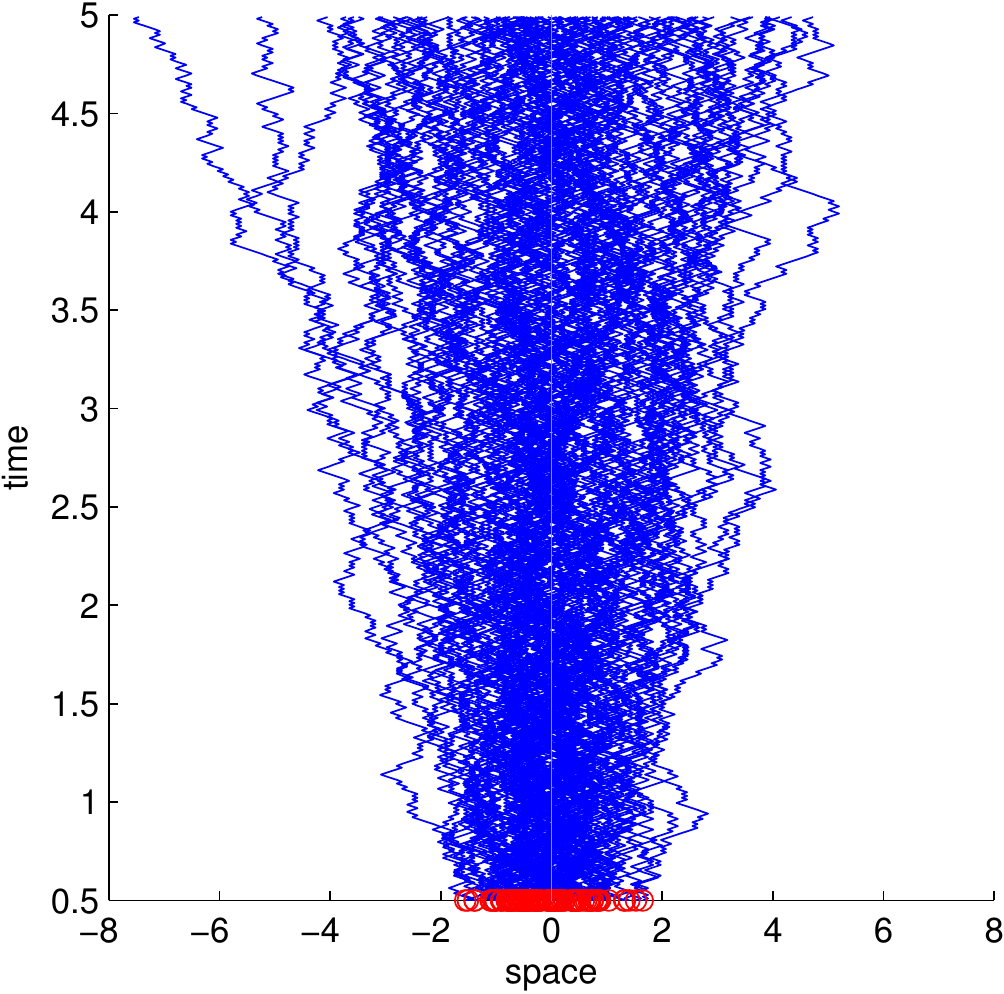}  & 
\includegraphics[width=0.23\textwidth]{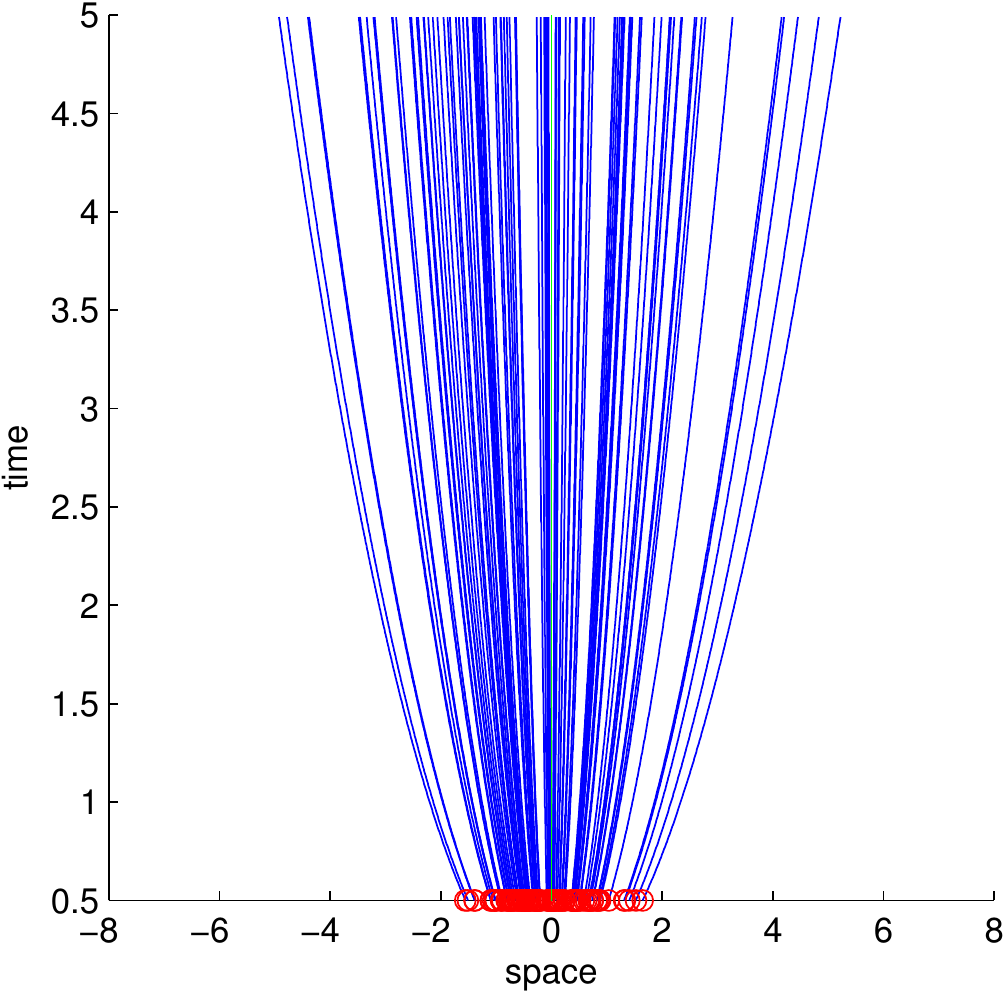} \\ 
(a) & (c)\\
\includegraphics[width=0.23\textwidth]{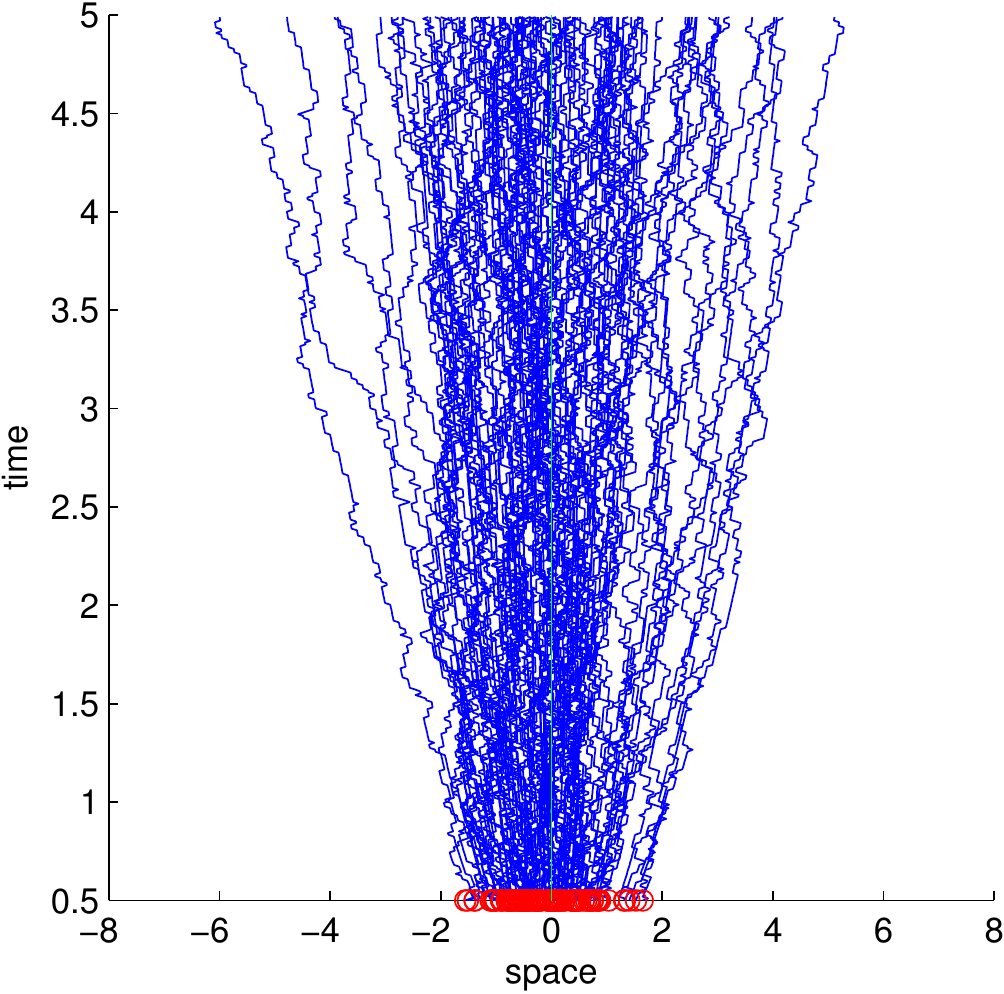} & 
\includegraphics[width=0.23\textwidth]{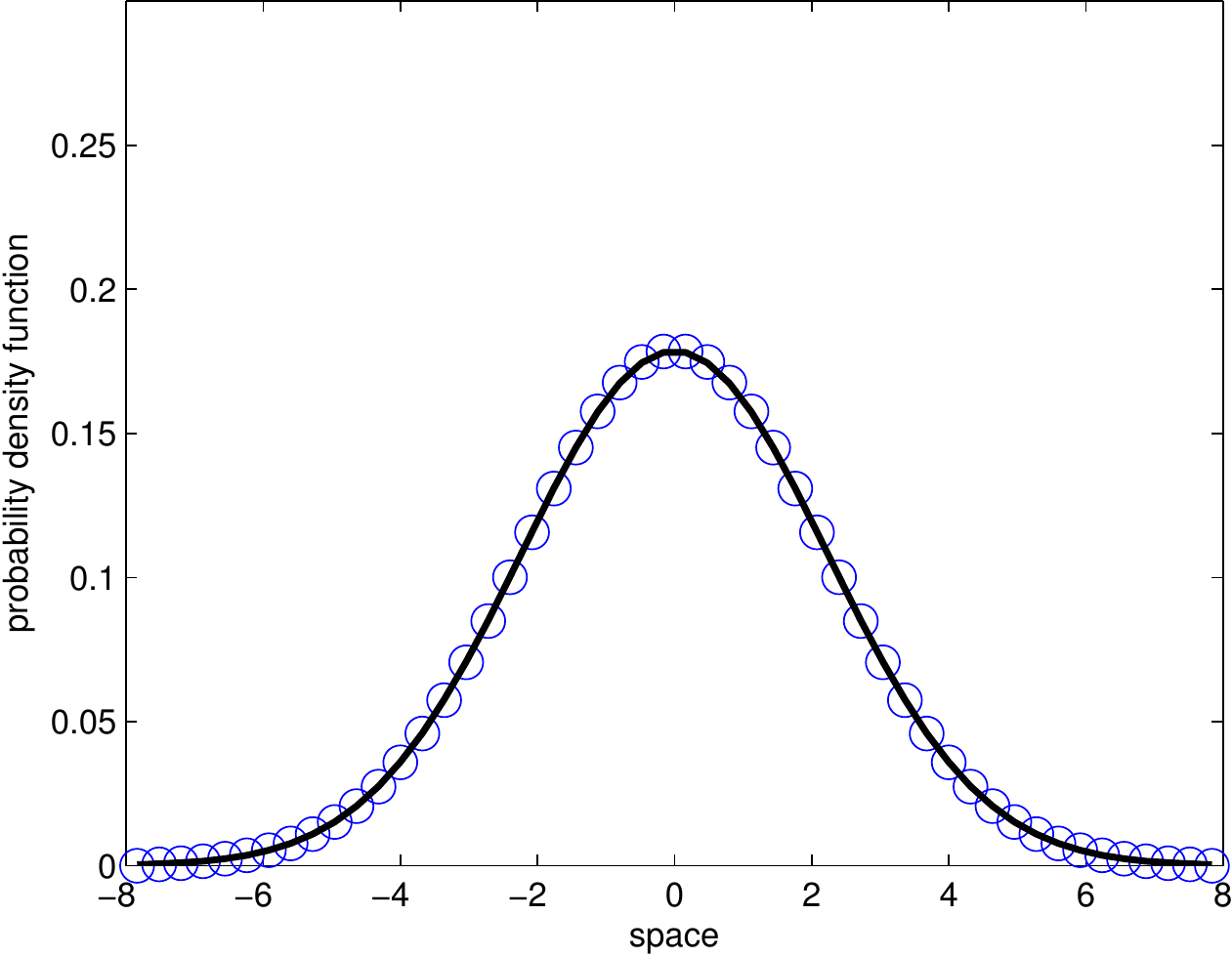} \\ 
(b) & (d)
\end{tabular}
\end{center}
\caption{Solution $\{Y^k\}_k$ to (\ref{eq:partialcouplingmicro}) in the plane $x$-$t$ for $N_p=100$ and (a) $\theta=1$, (b) $\theta=0.2$, (c) $\theta=0$. Red circles are the initial positions of the particles. (d) Probability density function $\{\psi^j_{N_t}\}_j$ of $N_p=500,000$ particles (blue circles) and the function $u(\tf,\cdot)$ (black solid line) for any $\theta$.}
\label{fig:partialcouplingmicro}
\end{figure}

Figure~\ref{fig:partialcouplingmicro} shows the solution to (\ref{eq:partialcouplingmicro}) (with $N_p=100$) and the corresponding probability density function $\psi$ (with $N_p=500,000$) at final time $\tf$ for $\theta=0,\ 0.2,\ 1$. We notice that the Brownian motion is regularized as expected since the trajectories are gradually more and more smooth as $\theta$ goes to 0, while numerical evidence shows that the overall statistical properties of the collective dynamics are kept.

\subsection{Partial coupling II: The Brownianized heat equation}\label{sec:partialcouplingmacro}
\begin{figure}[h!]
\begin{center}
\begin{tabular}{cc}
\includegraphics[width=0.23\textwidth]{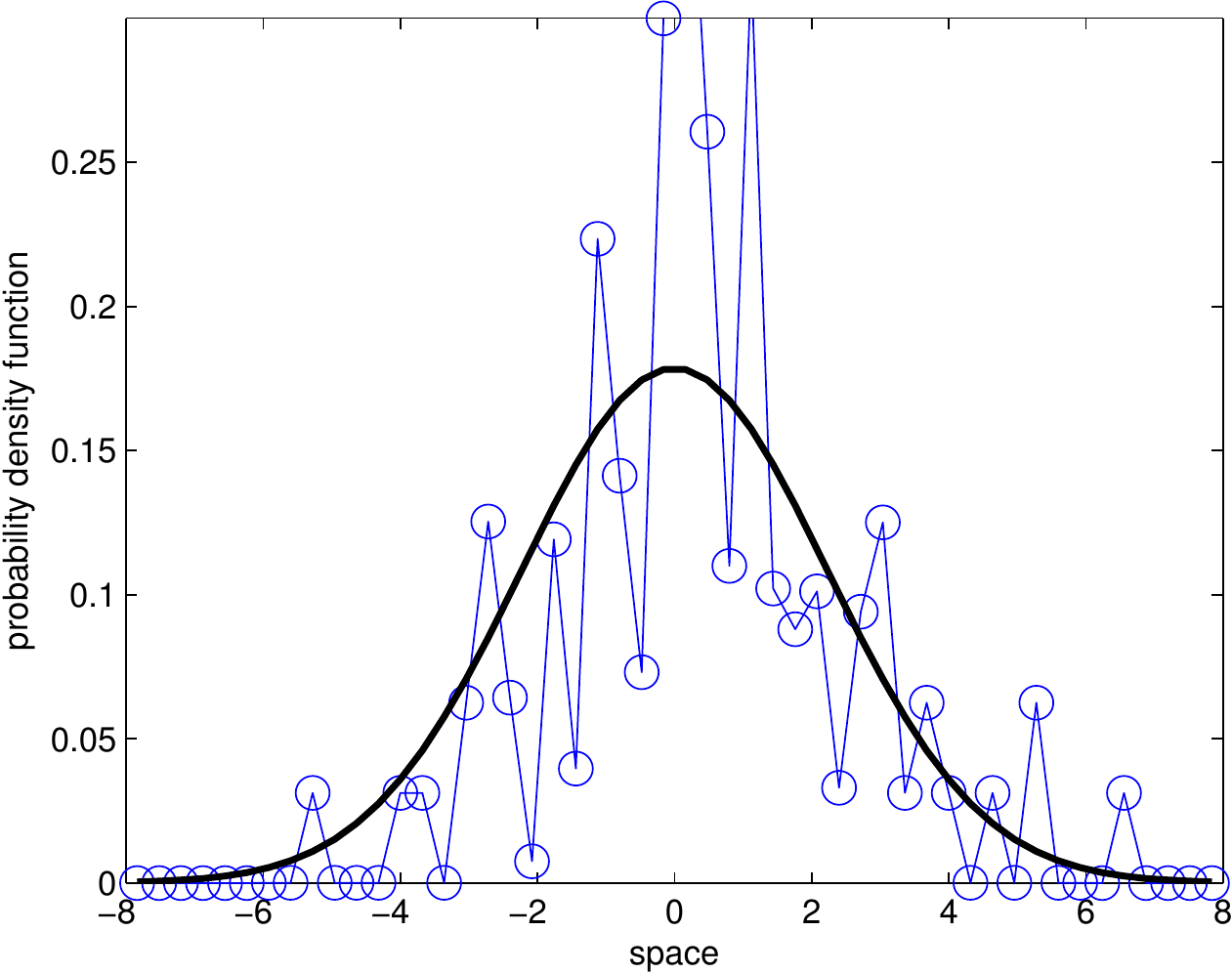}&
\includegraphics[width=0.23\textwidth]{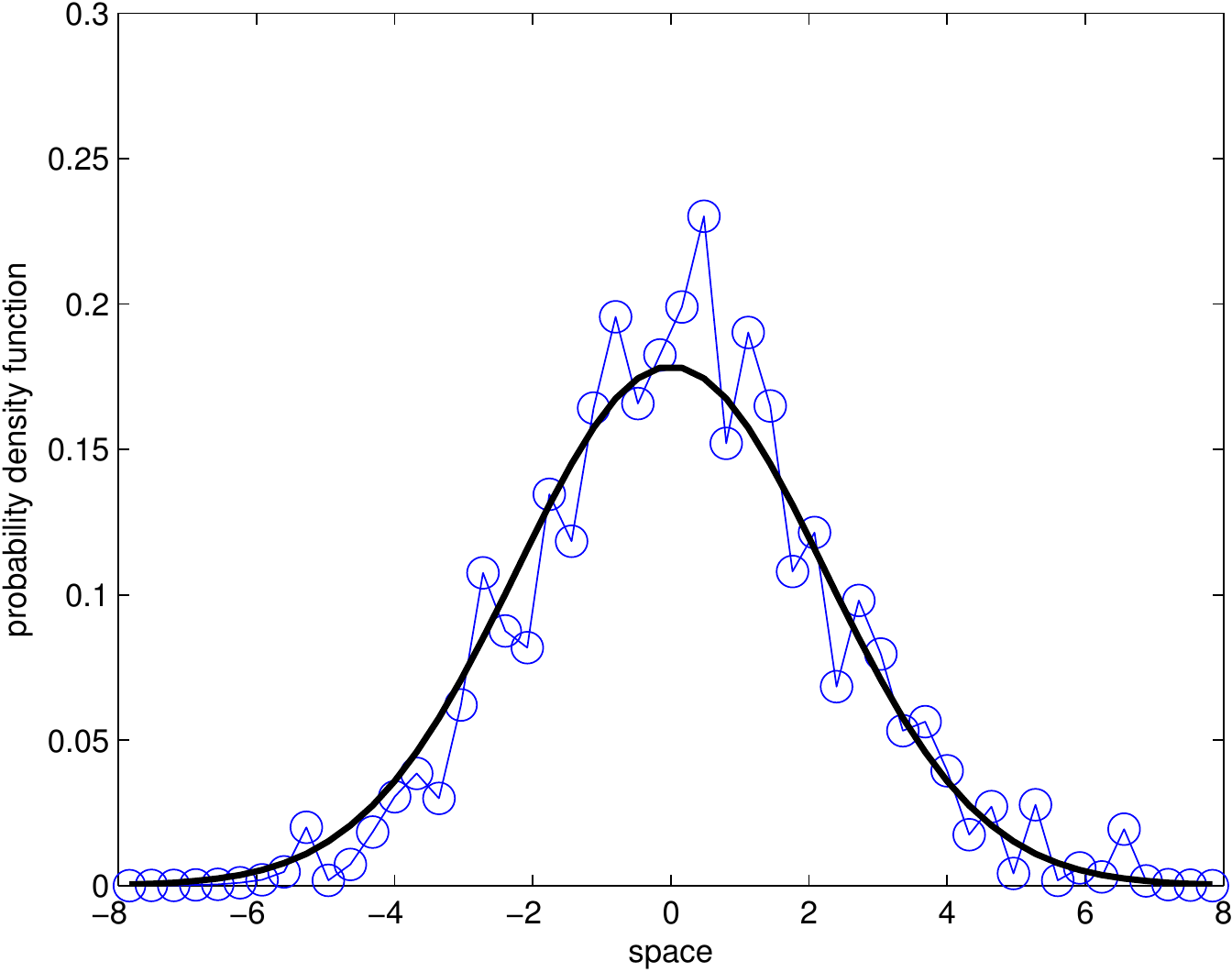}\\
(a) & (d) \\
\includegraphics[width=0.23\textwidth]{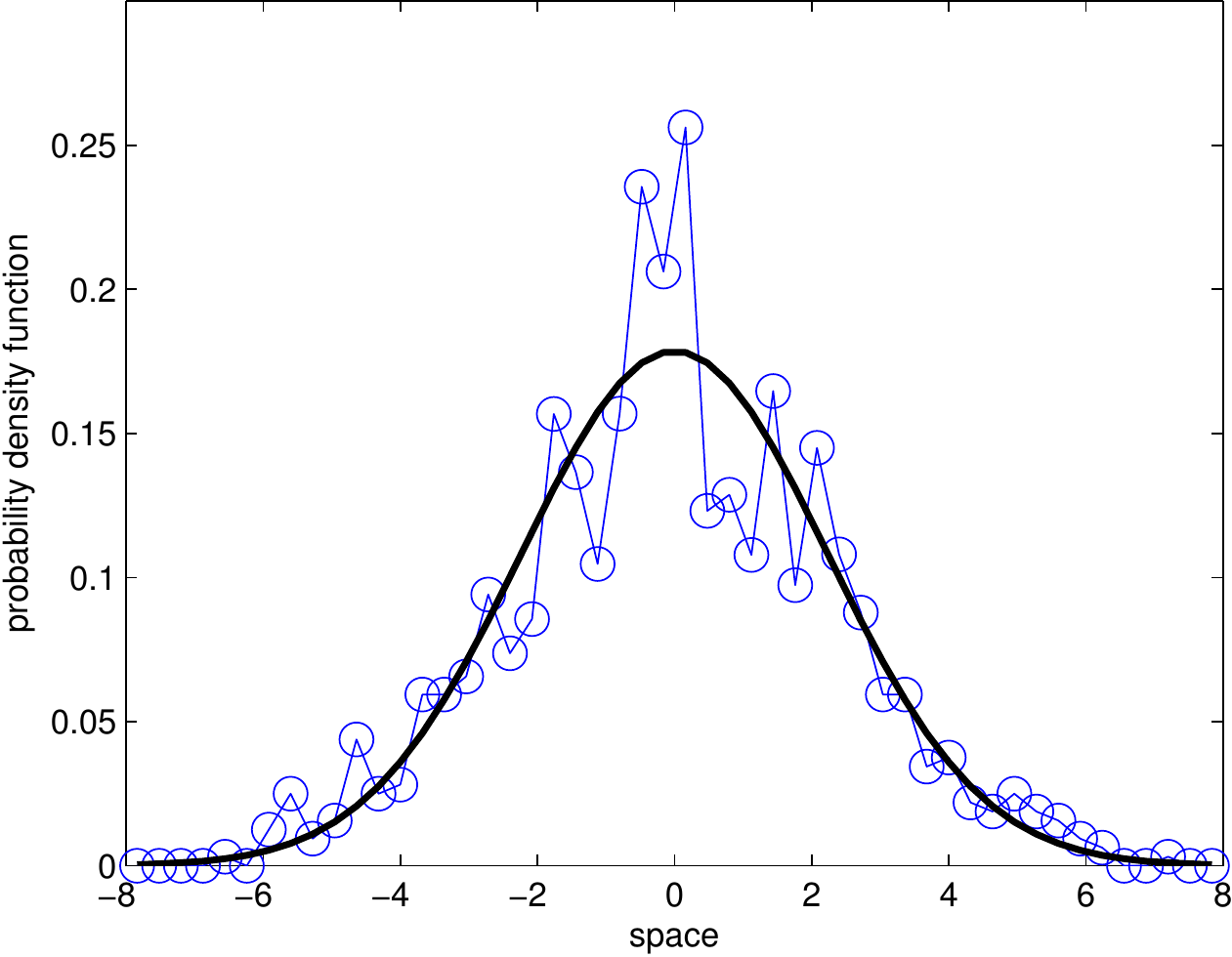}&
\includegraphics[width=0.23\textwidth]{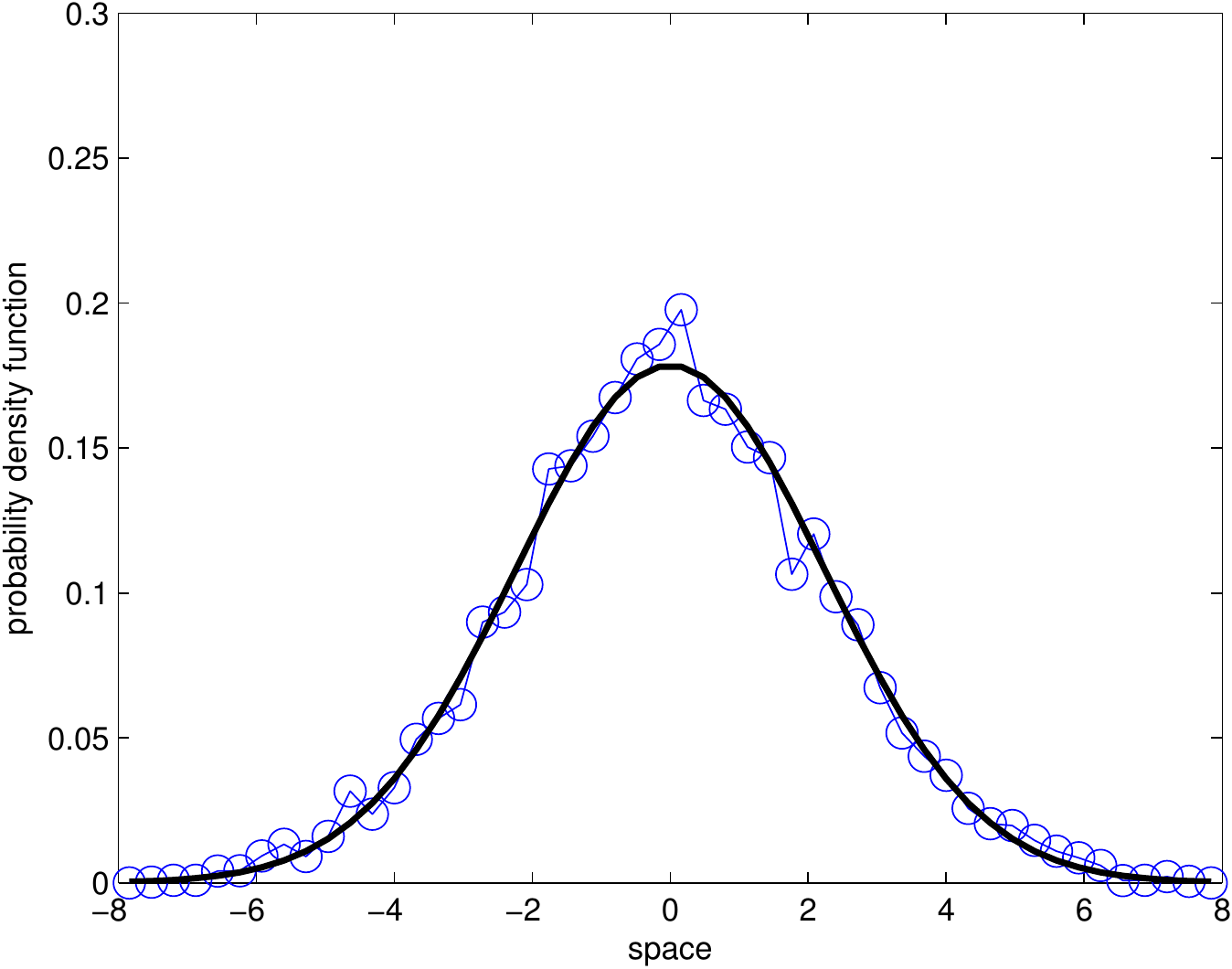}\\
(b) & (e) \\
\includegraphics[width=0.23\textwidth]{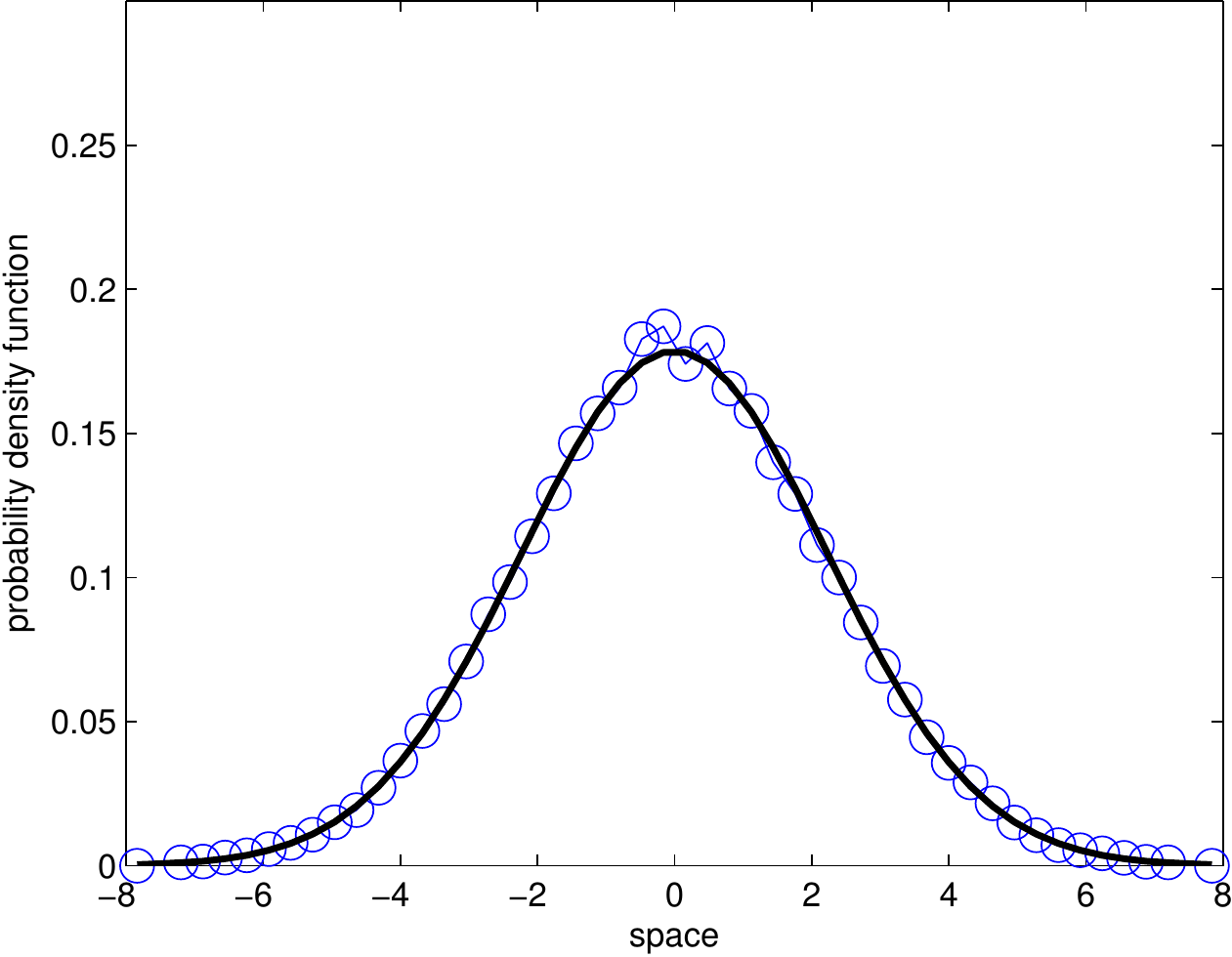} & 
\includegraphics[width=0.23\textwidth]{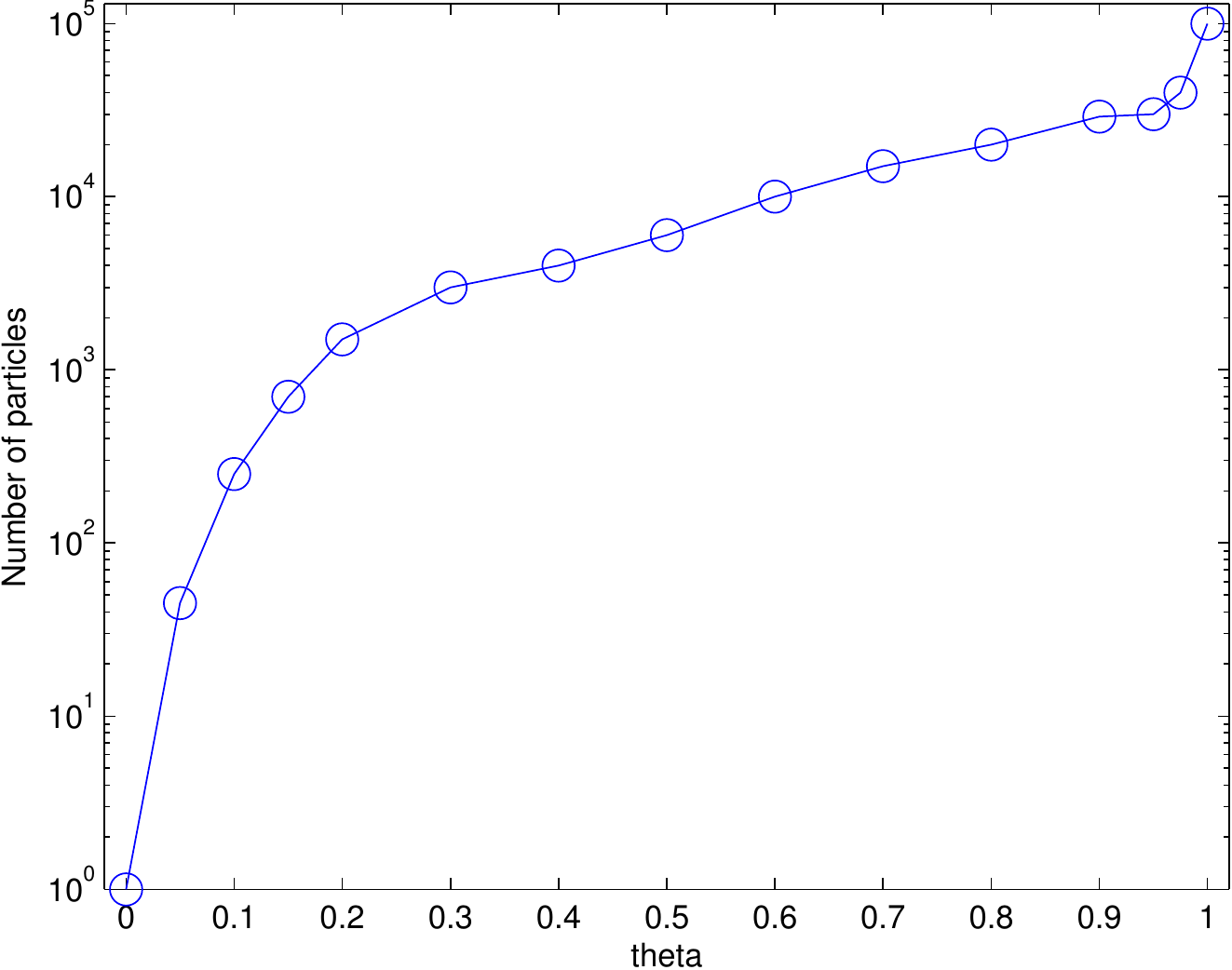} \\
(c) & (f)
\end{tabular}
\end{center}
\caption{Solution $\{\phi^j\}_j$ to (\ref{eq:schema_partialcoupling_macro}) at $t=\tf$ (blue circles) and the function $u(\tf,x)$ (black solid line) for (a) $\theta=1$ and $N_p=100$, (b) $\theta=1$ and $N_p=1000$, (c) $\theta=1$ and $N_p=100,000$, (d) $\theta=0.5$ and $N_p=100$, and (e) $\theta=0.5$ and $N_p=1000$. (f) $(\theta,N_p)$ pairs such that the solution $\phi_{N_t}$ to (\ref{eq:schema_partialcoupling_macro}) has $L^1$-distance from $u(\tf,\cdot)$ equal to 0.025, log scale on vertical axis.}
\label{fig:partialcouplingmacro}
\end{figure}

\noindent In order to obtain the opposite coupling, we first compute $\{Y_n^k\}_{n,k}$ by means of (\ref{Euler_for_SDEs}) and then we use the particles' trajectories to correct the evolution of $\phi$. To begin with, let us consider the extreme case $\theta=1$. Here the evolution of $\phi$ is totally driven by the particles. It is convenient to find a dynamics of the form 
$$
\phi_{n+1}^j=\phi_n^j+\Delta t S_n^j[Y_n,Y_{n+1}]
$$
for some function $S_n^j$ which depends on the particles' positions at time $t_n$ and $t_{n+1}$. At any fixed time step $n$, the function $S_n^j$ can be recovered \emph{a posteriori} as follows:
We define 
$$
\widetilde \phi^j_{n+1}[Y_n,Y_{n+1}]:=\phi_n^j+\sum_{k\in \mathcal I^j_n}\frac{1}{N_p \Delta x}-\sum_{k\in \mathcal O^j_n}\frac{1}{N_p \Delta x}
$$
where $\mathcal I^j_n$ is the set of particles which enter the cell $E^j$ in the time interval $[t_n,t_{n+1}]$ and $\mathcal O^j_n$ is the set of particles which leave the cell $E^j$ in the same time interval. In this way we account for the density which passes from one cell to another. Finally, we set
\begin{equation}\label{eq:S}
S_n^j[Y_n,Y_{n+1}]:=\frac{\widetilde \phi^j_{n+1}[Y_n,Y_{n+1}]-\phi_n^j}{\Delta t}.
\end{equation}
At this point the coupling follows easily by setting
\begin{multline}\label{eq:schema_partialcoupling_macro}
\phi_{n+1}^{j}=\phi_{n}^{j}+\Delta t\Big(
\theta S_n^j[Y_n,Y_{n+1}]+\\
(1-\theta)\frac{D}{\h^2}
\left(\phi_{n}^{j+1}-2\phi_{n}^{j}+\phi_{n}^{j-1}\right)\Big)
\end{multline}
for any $\theta\in[0,1]$. 
Figure \ref{fig:partialcouplingmacro} shows the solution $\phi$ to (\ref{eq:schema_partialcoupling_macro}) at final time $\tf$ for some values of $\theta$ and $N_p$. The comparison between Figures \ref{fig:partialcouplingmacro}(b) and \ref{fig:partialcouplingmacro}(d) suggests that \textit{the multiscale model can be effectively used to reduce the number of particles (and thus the numerical complexity) still keeping the same accuracy in the final solution}. Indeed, the solutions have approximately the same quality but they are obtained with $N_p=1000$ and $N_p=100$ particles respectively. This means that the macroscopic counterpart is able to compensate the lack of particles. 
To confirm this insight, in Figure \ref{fig:partialcouplingmacro}(f) we plot some pairs $(\theta,N_p)$ associated to the same (approximate) $L^1$-distance from $u$ at final time, more precisely we impose
$
E^1:=\sum_{j=1}^{N_x}|u(\tf,x^j)-\phi_{N_t}^j|\Delta x=0.025
$. 
Note that we use a log scale on vertical axis. 

\textit{Remark.} When the heat equation is solved by using microscopic information, it gets an hyperbolic flavor because the density is, at least partially, moved by means of the microscopic velocity field. Swapping the actors, this reminds the classical numerical schemes for the advection equation where an artificial viscosity is added. 
It is also interesting to note that, even if no CFL condition is required here, it can still be imposed. The condition is satisfied if particles remain in the same cell or move at most one cell apart in one time step $\Delta t$, i.e.\ $\Delta x \geq |\omega_n\sqrt{2D}\sqrt{\Dt}|=\sqrt{2D}\sqrt{\Delta t}$, which corresponds to the parabolic-type CFL condition for the scheme (\ref{schemeDFC}).

\subsection{Full coupling}\label{sec:fullcoupling}
\noindent The full coupling of Brownian motion and heat equation is obtained easily merging the ideas discussed in the previous sections. The main difference is that we have to approximate $v_M$ at each time step because it is no longer given by  (\ref{Vm_esplicito}). This can be done via a finite difference scheme, as
$$
v_M(t_n,Y_n^k)=-\frac{D}{u_n^{j^*}} \frac{(u_n^{j^*+1}-u_n^{j^*-1})}{2\Delta x}
$$
where $j^*=j^*(n,k)$ is defined as the unique index such that $Y_n^k\in E^{j^*}$. 
Summarizing, the scheme reads as
\begin{equation}\label{eq:fullcoupling}
\left\{
\begin{array}{l}
Y^k_{n+1}=
\left\{
\begin{array}{ll}
Y^k_n+\omega_n\sqrt{2D}\sqrt{\Dt}, & \textrm{with prob.\ } \theta, \\ [1mm]
Y^k_n-\frac{D\Delta t}{\phi_n^{j^*}} \frac{(\phi_n^{j^*+1}-\phi_n^{j^*-1})}{2\Delta x}, & \textrm{with prob.\ } (1\!-\!\theta),
\end{array}
\right.
\\ [6mm]
\phi_{n+1}^{j}=\phi_{n}^{j}+\Delta t\Big(\theta S_n^j[Y_n,Y_{n+1}]+\\ 
\hskip2.5cm (1-\theta)\frac{D}{\h^2}\left(\phi_{n}^{j+1}-2\phi_{n}^{j}+\phi_{n}^{j-1}\right)\Big)
\end{array}
\right.
\end{equation}
where $S_n^j$ is defined as in (\ref{eq:S}).

\begin{figure}[t!]
\begin{center}
\begin{tabular}{cc}
\includegraphics[width=0.23\textwidth]{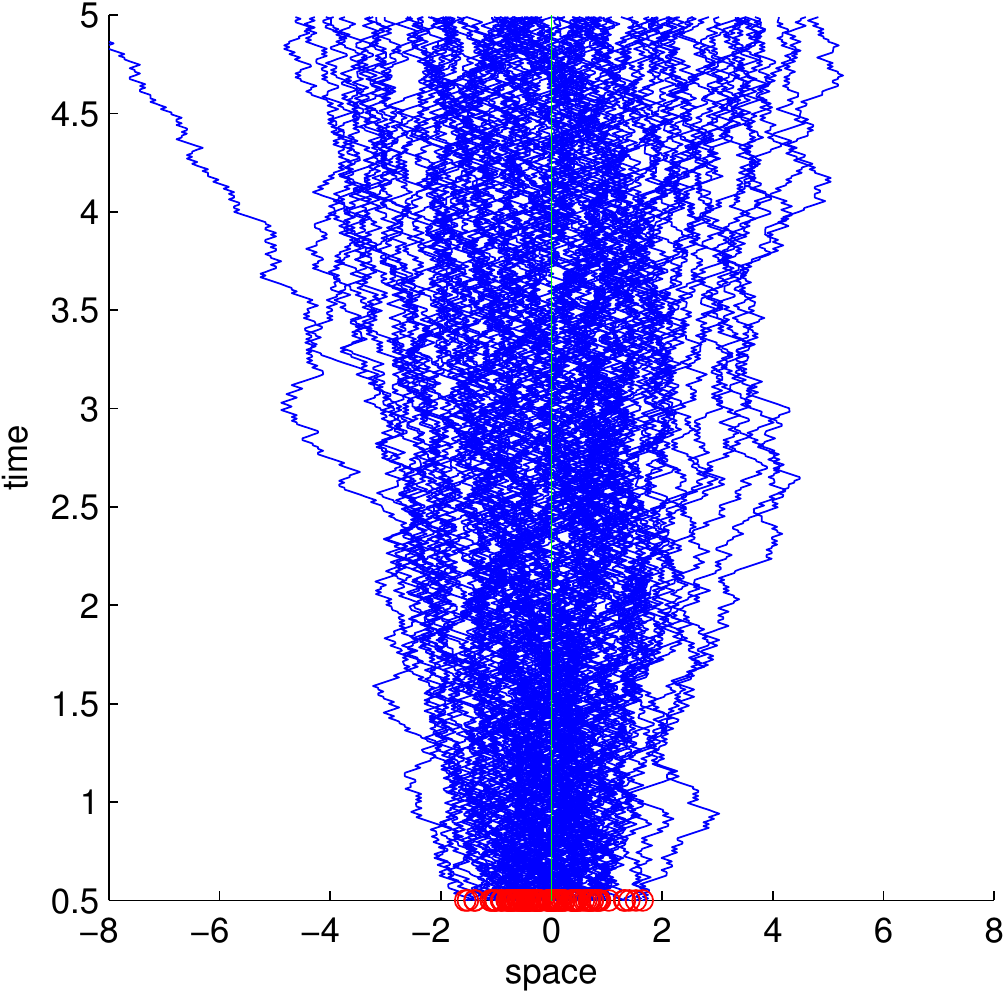} &
\includegraphics[width=0.23\textwidth]{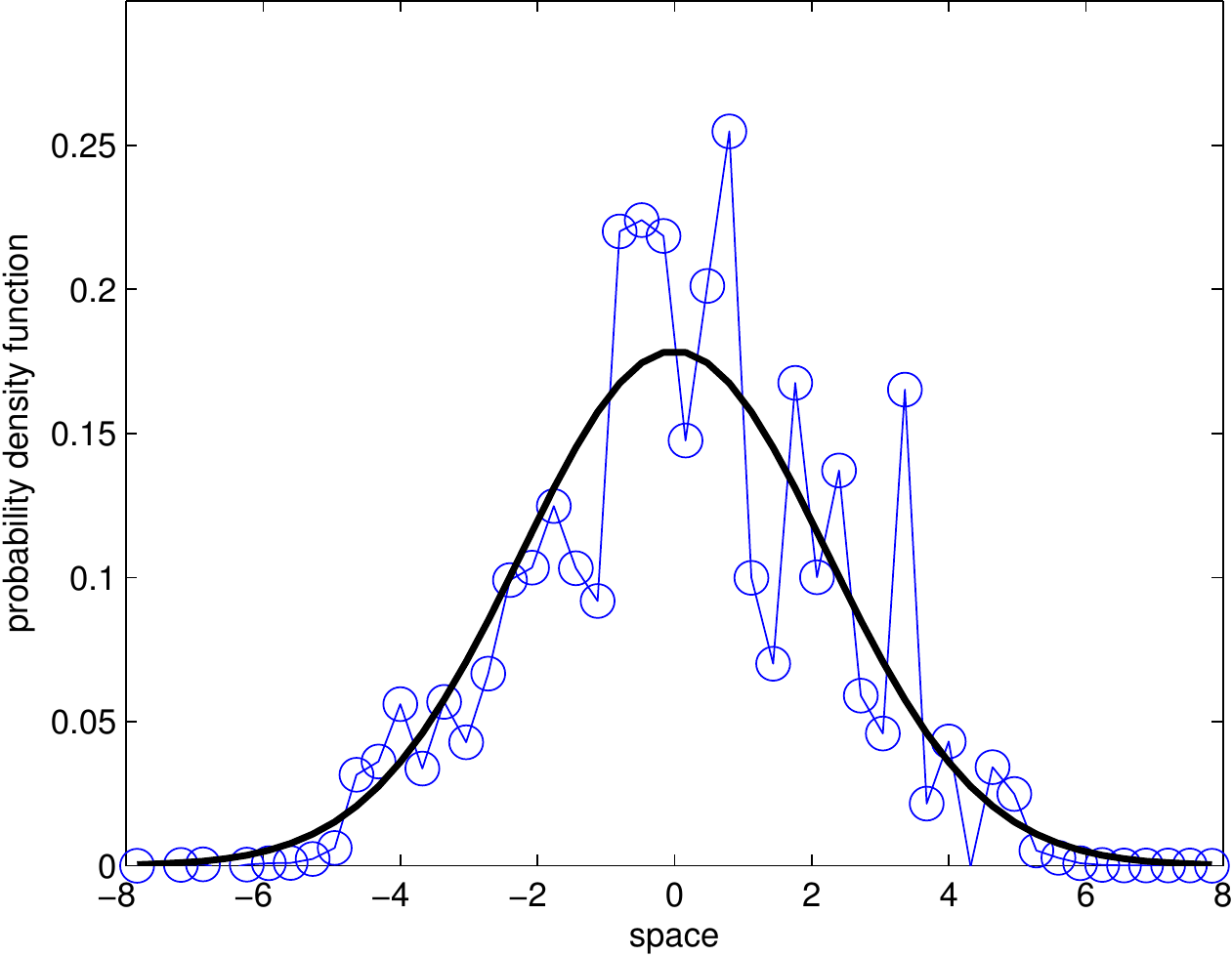}\\
(a)&(d)\\
\includegraphics[width=0.23\textwidth]{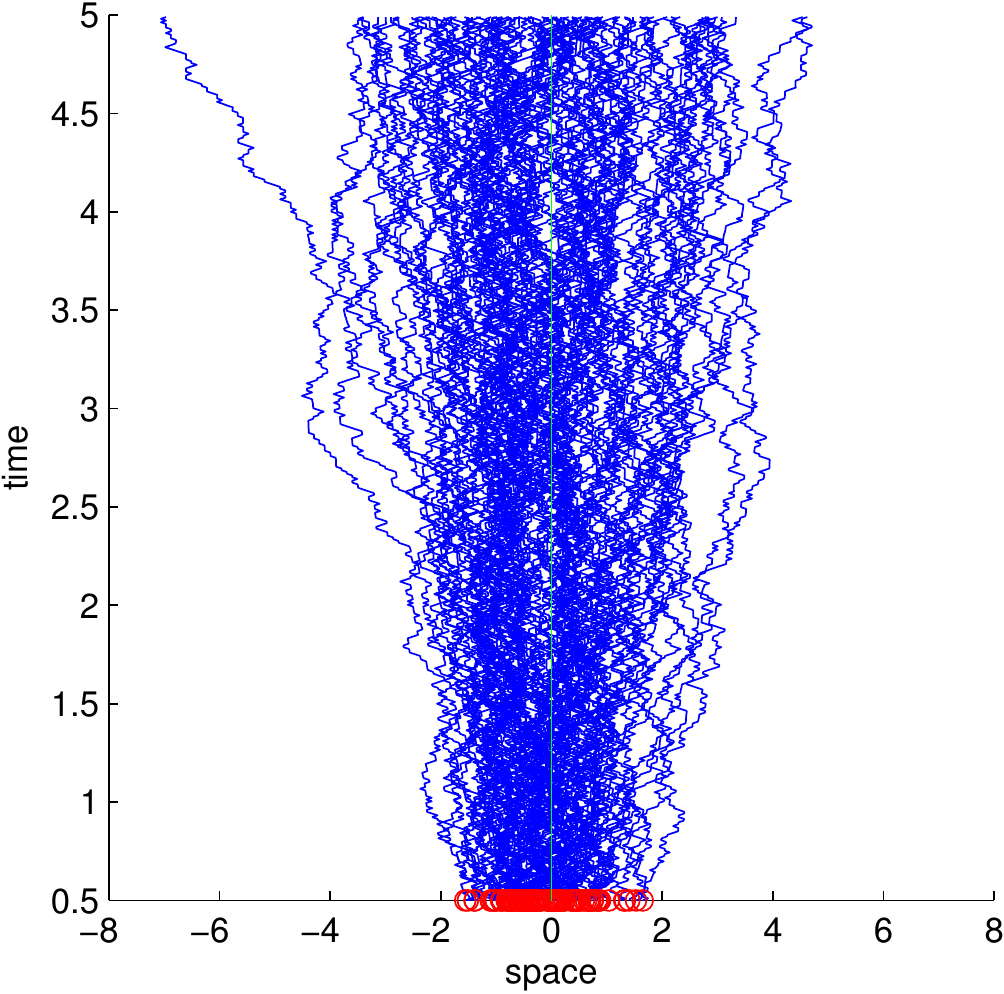} &
\includegraphics[width=0.23\textwidth]{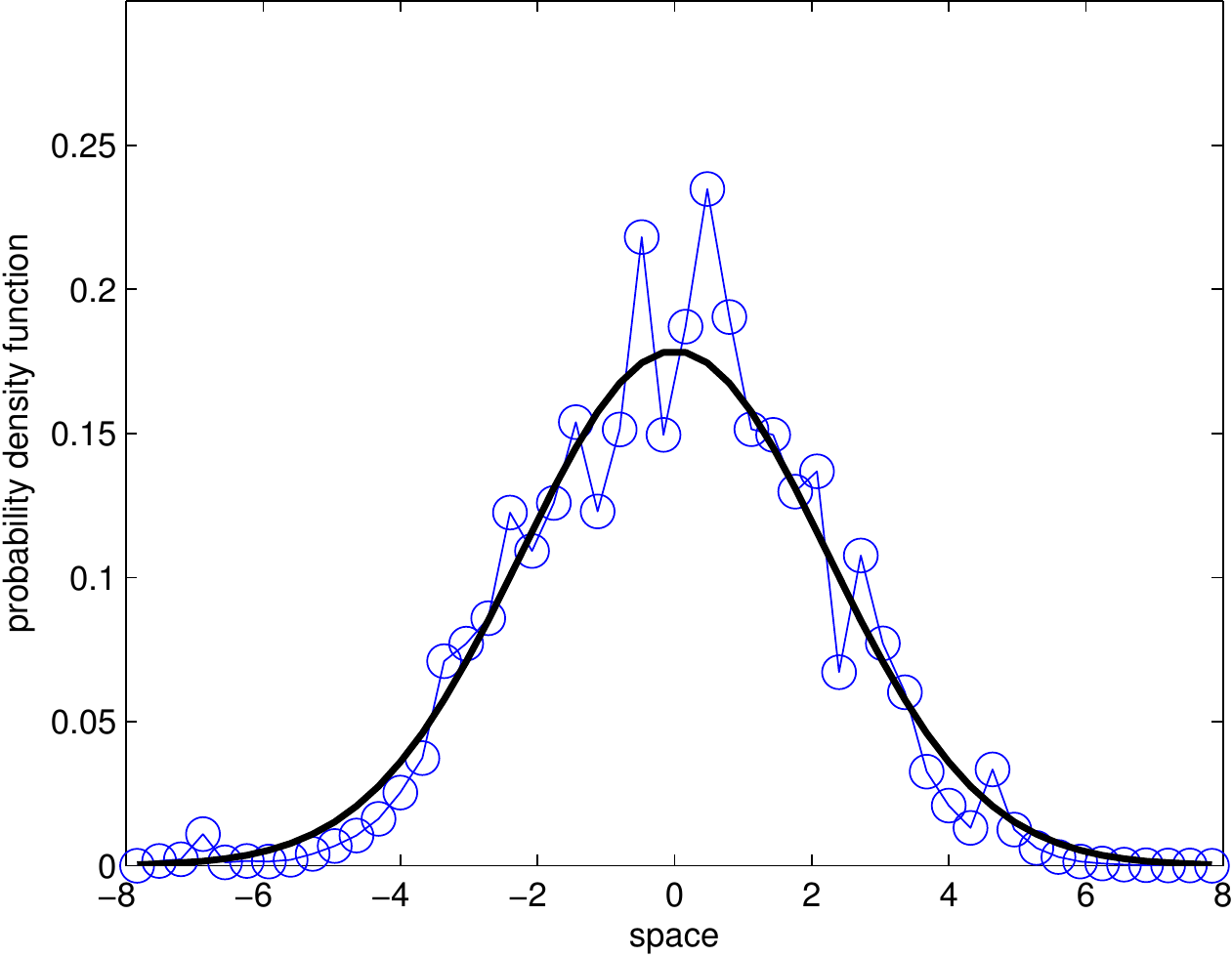}\\
(b)&(e)\\
\includegraphics[width=0.23\textwidth]{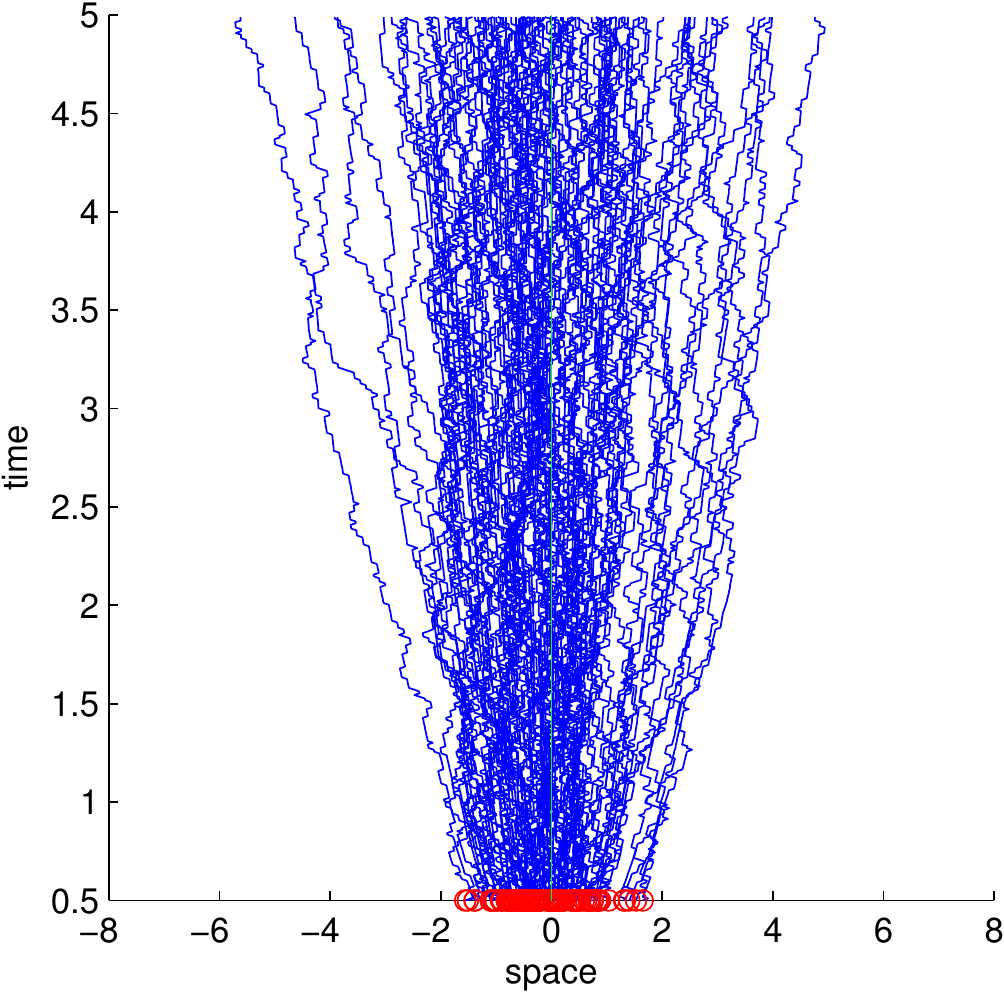} &
\includegraphics[width=0.23\textwidth]{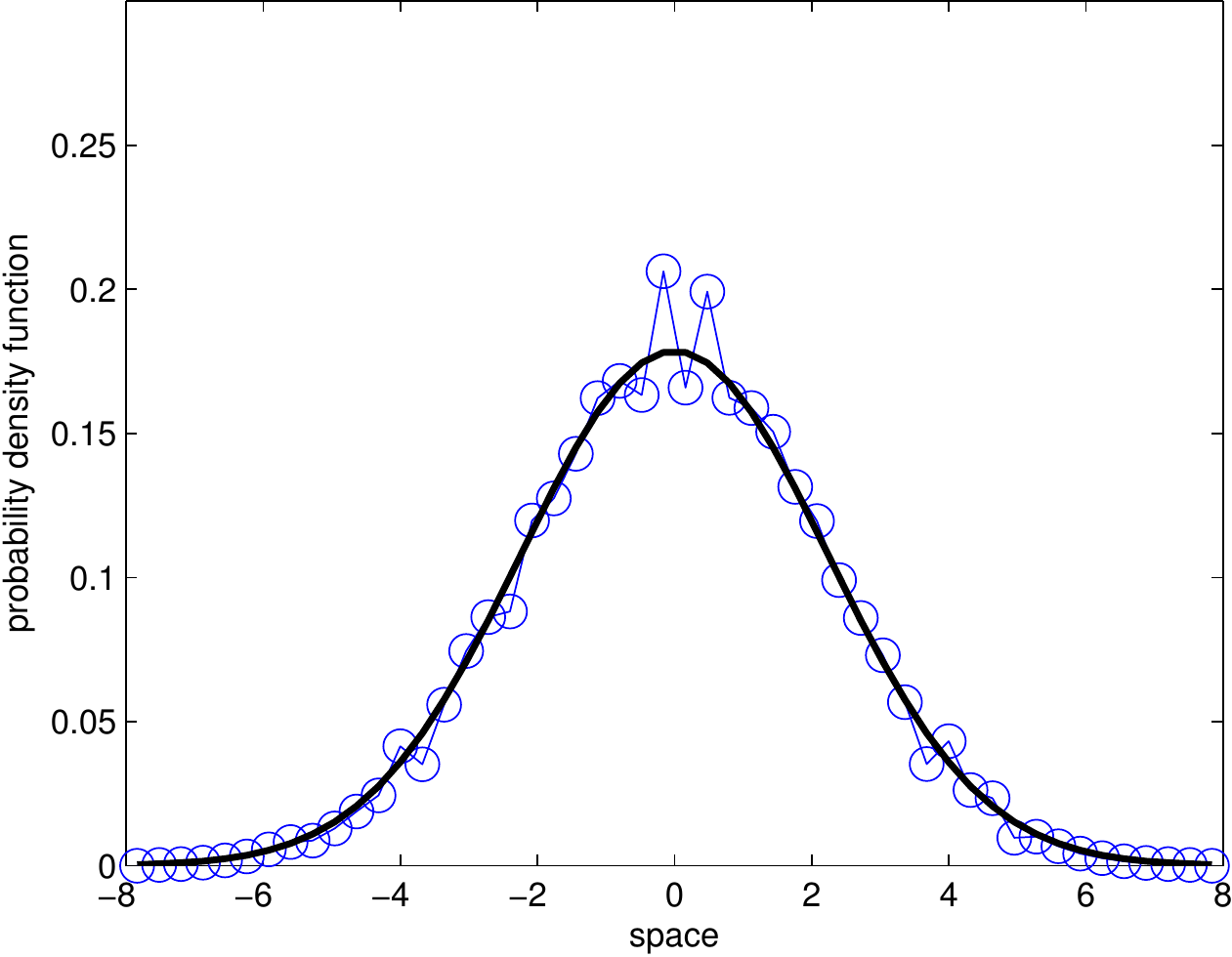}\\
(c)&(f)
\end{tabular}
\end{center}
\caption{Left column: Solution $\{Y^k\}_k$ to (\ref{eq:fullcoupling}) in the plane $x$-$t$ for $N_p=100$ and (a) $\theta=0.8$, (b) $\theta=0.5$, and (c) $\theta=0.2$.
Right column: Solution $\{\phi^j\}_j$ to (\ref{eq:fullcoupling}) at $t=\tf$ (blue circles) and the function $u(\tf,x)$ (black solid line), for $N_p=100$ and (d) $\theta=0.8$, (e) $\theta=0.5$, and (f) $\theta=0.2$.}
\label{fig:fullcoupling}
\end{figure}

Figure \ref{fig:fullcoupling} shows the solution $(Y,\phi)$ to (\ref{eq:fullcoupling}) for $N_p=100$ and $\theta=0.2$, $0.5$, and $0.8$. Here a genuine coupling is vi\-sible, while the global statistical properties are kept at any scale.

\section*{Final comments}
\noindent The coupling technique presented in this paper allows to blend the Brownian motion and the heat equation. The question arises why one should find convenient to blend the scales in this way. In \cite{cristiani2011MMS} the authors deal with granular flow and the power of the ``convex combination'' ruled by the parameter $\theta$ is evident, because the macrosolver is not able, alone, to catch the natural break of symmetry of the particles' density actually observed in reality, which is ultimately triggered by microscopic effects. If, instead, the macrosolver is able to describe the dynamics in full detail, as it happens here, the coupling seems not to represent a real value added. 
However, we think that the ideas presented here can be extended to more general diffusion-based equations modeling more complex phenomena, where the description at different scales \textit{is not perfectly symmetrical.} This means that the microscopic description is actually richer than the macroscopic one and single particles are responsible for the onset of some phenomenon (e.g., break of symmetry, fractures) which is then visible at large scale. 
In such new contexts, the constant presence of the microscale is crucial, but at the same time one wants to avoid tracking \textit{all} the particles, and prefers letting the macrosolver take the place of a part of them. The coupling seems also to be conceptually suitable to describe physical duality in quantum mechanics, even in the case of perfectly symmetric micro and macro descriptions. To do so, some ideas can be borrowed by the Lattice Boltzmann method \cite{succichapter,succi1993PhysD}.

A main open problem is given by the choice of $\theta$. Is there an optimal value for the coupling parameter? We think that the answer is necessarily problem-dependent. Likely, a good parameter $\theta=\theta(t,x)$ will be time- and space-dependent (it could the solution of an additional PDE). A first investigation in this direction can be found in \cite{colombi}. 

Let us also mention that our multiscale approach avoids to deal with classical micro-macro interfaces. The dynamics can be described at macrolevel ($\theta=0$) in a subomain and at macromicrolevel ($\theta\in(0,1]$) in another subdomain. Then, the hand-shaking at the boundary of the subdomains is done at macroscale level only (i.e.\ macro-with-macro, neglecting the microscale). In this way we get an easy seamless exchange of information, cf.\ \cite{e2009JCP}.

\section*{\normalsize Acknowledgments}
\label{Acknowledgements}
\noindent The author wants to thank S. Succi, B. Piccoli, A. Tosin, S. Cacace, and the anonymous referee for the useful suggestions, and A. Einstein for the inspiration.



\end{document}